\documentclass[12pt]{amsart}
\usepackage{amsfonts,amssymb}
\usepackage{graphicx}
\usepackage{amscd}
\usepackage{eucal}
\input epsf
\def\Z{\mathbb{Z}}
\def\C{\mathbf{C}}

\def\PSL{{\rm PSL}}
\def\<{\langle }
\def\>{\rangle }

\begin{document}
\setlength{\baselineskip}{18pt}
\title[]
{Cesaro averages of  Euler-like  functions }

\author[]{F. Aicardi}
%
%
%

\begin{abstract}
   By  Euler-like  function  we mean  a function defined on the positive integers and associating  to  $n$ the product, over all primes  $p$ dividing $n$, of   1  plus
(or minus)  the inverse of $p$  to the power $s$.   We   calculate  the limit
of  the Cesaro means of these functions.
\end{abstract}
\maketitle

\section{Origin of the problem}

When work  [1] was in progress,  the author asked me about the Cesaro mean of  the  following rational  function $f$ of the natural $n$ \[ f(n)= \prod_{p|n}\left(1+\frac{1}{p}\right), \quad p \ \hbox{prime}, \]
(the integer $nf(n)$ being  the  number  of  straight   lines   through the origin in the torus  $\Z_n^2$, permuted  by   the action of $\PSL(2,\Z_n)$).
For  $n$ till  $10^7$,  I  found  that   the   Cesaro  mean  of $f$,  equal to  $ 1.5198177542107...$, multiplied by $\pi^2$,  becomes
$14.9999999958...$.   To prove that  this Cesaro  average  converges indeed to  $15/\pi^2$, I studied the Cesaro means of  more general functions,  that I call `Euler-like', since the  Euler  function $\phi(n)$ is equal to $n\prod_{p|n}\left(1-\frac{1}{p}\right)$ ($p$ prime).   I  hope these results   can be applied to generalize  some Arnold's theorems  in [1]. The  idea of   generalizing  the present  results to       Dirichlet L-series was suggested  by  V.Timorin.

\section{Two theorems}
{ \bf Theorem 1.} {\it For every  natural\footnote{In fact, this  statement holds  for any complex number $n$ with real part greater than  1.} $n> 1$, the limit for $N\rightarrow \infty$ of the following Cesaro mean :
\begin{equation}\label{zetaplus}  \frac{1}{N} \sum_{m=1}^N \prod_{p|m} \left(1+\frac{1}{p^{n-1}}\right)
\end{equation}
 is equal to
 \[  \zeta(n)/\zeta(2n) .\] }

 In particular, for $n=2$,
 \begin{equation}\label{quindici}  \lim_{N\rightarrow \infty} \frac{1}{N} \sum_{m=1}^N \prod_{p|m} \left(1+\frac{1}{p}\right)
= \frac{\zeta(2)}{\zeta(4)} = \frac{\pi^2/6}{\pi^4/ 90}= \frac {15}{\pi^2}. \end{equation}

{ \it Proof} Remember that
\begin{equation}\label{zetan}    \prod_p \left(1-\frac{1}{p^n}\right)= \frac{1}{\zeta(n)},
\end{equation}
where the product is extended to all primes.

Hence (see also \cite{Har}, \cite{Mon})
\begin{equation} \label{zetan2}    \prod_p \left( 1+\frac{1}{p^n}\right) =
\frac{ \prod_p \left(1-\frac{1}{p^{2n}}\right)}{\prod_p \left(1-\frac{1}{p^n}\right)}=\frac {\zeta(n)}{\zeta(2n)}
 \end{equation}
where the products are  extended to all primes.

We have
\[ \lim_{N\rightarrow \infty} \frac{1}{N} \sum_{m=1}^N \prod_{p|m} \left(1 +\frac{1}{p^{n-1}}\right)=
  \lim_{N\rightarrow \infty}   \frac{1}{N} \sum_{m=1}^N \left(\sum_{k_\circ|m} \frac {1}{k_\circ^{n-1}}\right) \]
where $k_\circ$  is  a natural  {\it square  free}  (i.e., a natural   number  which is not divisible by the square of any prime number), and the    $\sum_{k_\circ|m} \frac {1}{{k_\circ}^{n-1}}$ is extended  to all $k_{\circ}$ dividing $m$.

The number of values of $m\le N$  divisible by $k$ is  $[N/k]$, approaching
$N/k$ for    for $ N\rightarrow \infty$. Hence:
\[ \lim_{N\rightarrow \infty}   \frac{1}{N} \sum_{m=1}^N \left(  \sum_{k_\circ|m} \frac {1}{k_\circ^{n-1}}\right)= \lim_{N\rightarrow \infty}   \sum_{k_\circ|m} \frac{1}{N} \frac {[N/k_\circ]} {k_\circ^{n-1}}=    \sum_{k_\circ} \frac {1}{k_\circ^{n}},\]

We write finally
\[ \sum_{k_\circ} \frac {1}{k_\circ^{n}}= \prod_p \left(1+\frac{1}{p^n}\right),  \]
since this last infinite product is the sum of all numbers which are powers of simple products  of primes,
as the sum over all $k_\circ$   at the first  member  of equality.
The  claim of the  theorem is obtained using  eq.  (\ref{zetan2}).

We state now a similar theorem,  with two proofs. The first proof  holds only for naturals $n>1$, while
the second one holds for  any $n$ with real part greater than 1.

{\bf  Theorem 2.} {\it For every  natural $n$, the limit for $N\rightarrow \infty$ of the following Cesaro mean is equal to the inverse of the $\zeta(n)$:}
\begin{equation}\label{zetasomma}  \lim_{N\rightarrow \infty} \frac{1}{N} \sum_{m=1}^N \prod_{p|m} \left(1-\frac{1}{p^{n-1}}\right)=\frac{1}{\zeta(n)}.
\end{equation}

{\it First proof.}
Formula at the left member of (\ref{zetan})  can be  interpreted  as  the probability   that  $n$  numbers $m_1,m_2,\dots, m_n$ randomly  chosen  in  $\mathbb{N}$
(with uniform probability distribution)  satisfy   $\gcd(m_1,m_2,\dots, m_n)=1$.

Indeed, the probability   that  a natural  number  $m$  be divisible by $p$ is
equal to $ \frac{1}{p}$.

The probability that $n$ natural numbers  be all divisible by  $p$ equals  $(1/p)^{n}$, since
the  events (for different $m$)  ``{\it $m$ is divisible by the prime $p$}" are independent.

The probability that, for a given $p$, the $n$ chosen numbers be not all divisible by $p$ is  $(1-1/p^n)$.   The probability  that the $n$ chosen
numbers be not all divisible for any prime less or equal to $N$ is   $\prod_{p} \left(1-\frac{1}{p^n}\right)$, since the  events (for different primes $p$)  ``{\it  the $n$ chosen naturals are not all divisible by $p$}"  are independent.  The last expression is  the  left member of   (\ref{zetan}).

We  prove that  the left member of (\ref{zetasomma}) expresses the same probability as the left member of
(\ref{zetan}).
Suppose  now that $n$  numbers be  arbitrarily chosen  in the interval $I_N=[1,2,\dots,N]$ with uniform probability distribution.   Let $m$ be the maximum of  them.
The  probability  of choosing $m$,  as the probability of choosing any other number,  is equal to $1/N$.
The probability  that  the $n-1$ other  numbers be not all divisible by a prime $p$  dividing $m$ equals
  \[  \prod_{p|m} \left(1-q^{n-1}\right), \]
where  $q$   is  the probability   that the number  $m\in I_N$  be divisible by $p$:
\[   q:= p_{m\le N}(p|m)= \left\lfloor \frac{N}{p} \right \rfloor /N. \]
  The probability  that the  $n$  chosen  numbers in $I_N$ be not all divisible by a common prime is the
  sum, over all the values    $m$ of the maximum among these numbers,   of the  probabilities that the other $n-1$ numbers be not all  divisible by the  primes dividing $m$, i.e.:
 \begin{equation} \label{pro} \frac{1}{N} \sum_{m=1}^N \prod_{p|m} \left(1-q^{n-1}\right). \end{equation}
 Since
 \[  \lim_{ N\rightarrow \infty} q=\frac{1}{p},  \]
 the limit for $N\rightarrow \infty$  of (\ref{pro}) is equal to  the probability expressed by the left member of (\ref{zetan}),  and we thus obtain     eq. (\ref{zetasomma}).  \hfill $\square$

{\it Second proof.}  We have
\[ \lim_{N\rightarrow \infty} \frac{1}{N} \sum_{m=1}^N \prod_{p|m} \left(1-\frac{1}{p^{n-1}}\right)=
  \lim_{N\rightarrow \infty}   \frac{1}{N} \sum_{m=1}^N \left(\sum_{k_\circ|m} \frac {(-1)^{f_k}}{k_\circ^{n-1}}\right) \]
 where  $k_\circ$ is  square  free  and  $f_k$ is the number of primes entering the factorization
of $k_\circ$.

By the same  arguments used in the proof of Theorem 1:
\[ \lim_{N\rightarrow \infty}   \frac{1}{N} \sum_{m=1}^N \left(  \sum_{k_\circ|m} \frac {(-1)^{f_k}}{k_\circ^{n-1}}\right)= \lim_{N\rightarrow \infty}   \sum_{k|m} \frac{1}{N} \frac {(-1)^{f_k}[N/k_\circ]} {k_\circ^{n-1}}=    \sum_{k_\circ} \frac {(-1)^{f_k}}{k_\circ^{n}},\]
 and  write finally
\[ \sum_{k_\circ} \frac {(-1)^{f_k}}{{k_\circ}^{n}}= \prod_p \left(1-\frac{1}{p^n}\right).  \]
 The  claim of the  theorem is obtained using  eq.  (\ref{zetan}).  \hfill $\square$

\section{Generalization to Dirichlet L-series}
Let  $\chi$  be any  Dirichlet  character\footnote{ A Dirichlet character is any function $\chi: \Z  \rightarrow \C$  such  that:  1) it is periodic:   there exists a positive integer $k$ such that $\chi(n) = \chi(n + k)$ for all $n$; 2) $\chi(n)=0$ iff  $\gcd(n,k)>1$; 3) it is  multiplicative:  $\chi(mn) = \chi(m)\chi(n)$ for all integers $m$ and $n$}    $s$ a complex number with real part  greater than 0, and
\[   L(\chi,s)= \sum_{m=1}^\infty  \frac{\chi(m)}{m^s}  \]
the  corresponding   Dirichlet  L-series.

We recall that  the square of  a character is a character  as well.

{\bf  Theorem 3.}  {\it  The  limit for $N\rightarrow \infty$ of the  following Cesaro mean:
 \begin{equation}   \label{Dirminus}  \frac{1}{N} \sum_{m=1}^N \prod_{p|m} \left(1-\frac{\chi(p)}{p^{s-1}}\right)
\end{equation}
 is equal to
 \[  \frac{1}{L(\chi,s)} .\]
 The  limit for $N\rightarrow \infty$ of the  following Cesaro mean:
 \begin{equation}   \label{Dirplus}  \frac{1}{N} \sum_{m=1}^N \prod_{p|m} \left(1+\frac{\chi(p)}{p^{s-1}}\right)
\end{equation}
 is equal to }
 \[  \frac{ L(\chi,s)}{L(\chi^2,2s)} .\]

{\it Proof.}

We have only to  remark that
\[    L(\chi,s)= \frac{1} {\prod_p \left( 1- \frac {\chi(p)}{p^s} \right)}, \]
since  any  character $\chi$ is multiplicative  and the  product is extended to all primes $p$.

Hence
\begin{equation} \prod_p \left( 1+ \frac {\chi(p)}{p^s} \right)= \frac{\prod_p \left( 1- \frac {\chi^2(p)}{p^{2s}} \right)}{\prod_p \left( 1- \frac {\chi(p)}{p^s} \right)}= \frac{L(\chi,s)}{L(\chi^2,2s)} .
\end{equation}

Moreover
\[ \prod_p \left( 1+ \frac {\chi(p)}{p^s} \right)=\sum_{k_\circ} \frac{\chi(k_\circ)}{{k_\circ}^s}, \]
where the last sum is extended to all values of $k_\circ$ square free.

The limits for $N\rightarrow \infty$  of (\ref{Dirminus}) and (\ref{Dirplus})  are  then obtained the same  way  as the second
proof of Theorem  2 and as the proof of Theorem 1, respectively.


\end{document}